\documentclass[a4paper,11pt]{article}

\usepackage[utf8]{inputenc}
\usepackage[normalem]{ulem}
\usepackage{amsmath}
\usepackage{amsthm}
\usepackage{amsfonts}
\usepackage{amssymb}
\usepackage{mathrsfs}
\usepackage{graphicx}
\usepackage{enumerate}
\usepackage{url}
\usepackage{hyperref}
\usepackage{color, soul}
\usepackage[usenames,dvipsnames]{xcolor}
\usepackage{authblk}
\usepackage{geometry}
\usepackage{tikz}
\usepackage{caption}
\usepackage{subcaption}
\usepackage{mathtools}
\usepackage[numbers,sort]{natbib}

\usetikzlibrary{positioning,calc}
\tikzset{
    vertex/.style={
        circle,fill=black,scale=0.5
        },
    blank/.style={
        circle,fill=white
        },
    }

\usetikzlibrary{arrows}
\tikzset{
    vertex/.style={
        circle,fill=black,scale=0.5
        },
    blank/.style={
        circle,fill=white
        },
    }

\newcommand{\dgr}[1]{\gamma_{gr}^{\times 2}\left(#1\right)}
\newcommand{\dgri}[1]{\gamma_{gr}^{\leq 2}\left(#1\right)}
\newcommand{\gr}[1]{\gamma_{gr}\left(#1\right)}
\newcommand{\g}[1]{\gamma\left(#1\right)}

\newcommand{\N}{\mathbb{N}}

\newtheorem{theo}{Theorem}[section]
\newtheorem{lemma}[theo]{Lemma}
\newtheorem{prop}[theo]{Proposition}
\newtheorem{coro}[theo]{Corollary}

\theoremstyle{definition}

\newtheorem{rem}[theo]{Remark}

\title{Grundy double domination number: bounds, graph operations, and efficient computation for $P_4$-tidy graphs\footnote{Partially supported by PIP CONICET 1900, PICT-2020-03032 and PID 80020210300068UR.} }

\author{Pablo Torres\footnote{e-mail: ptorres@fceia.unr.edu.ar}}

\date{ }

\affil{\footnotesize Facultad de Ciencias Exactas, Ing. y Agrimensura, Universidad Nacional de Rosario, and Consejo Nacional de Investigaciones Científicas y Técnicas, Argentina}

\begin{document}

\maketitle

\begin{abstract}
Inspired by graph domination games, various domination-type vertex sequences have been introduced, including the Grundy double dominating sequence (GDDS) of a graph and its associated parameter, the Grundy double domination number (GDDN). The decision version of the problem of computing the GDDN is known to be NP-complete, even when restricted to split graphs and bipartite graphs. In this paper, we establish general tight bounds for the GDDN. We also describe GDDSs for vertex-removed graphs and for the join of two graphs. Applying these results, we prove that computing the GDDN is linear for $P_4$-tidy graphs, thereby solving an open problem previously posed for cographs by B.~Bre\v{s}ar et al. in~\cite{BPS-2021}.
\end{abstract}

\section{Introduction}
All graphs considered in this paper are simple and undirected. Given a graph $G=(V(G),E(G))$, the \emph{open neighborhood} of a vertex $v$ is the set $N_G(v) = \{u \in V(G):\ vu\in E(G)\}$.
The \emph{degree} of a vertex $v$ in $G$, denoted $deg_G(v)$, is the number of neighbors, $|N_G(v)|$, of $v$ in $G$. 
The \emph{minimum degree} and maximum degree among all the vertices of $G$ are denoted by $\delta(G)$ and $\Delta(G)$, respectively.
The \emph{closed neighborhood} of a vertex $v$ is the set $N_G[v]=N_G(v)\cup\{v\}$. A \emph{pendant} (resp. \emph{isolated}) vertex of a graph is a vertex of degree~$1$ (resp. degree~$0$). We denote by $I(G)$ the set of isolated vertices of a graph $G$. Two vertices $u,v\in V(G)$ such that $N_G[u]=N_G[v]$ (resp. $N_G(u)=N
_G(v)$) are called {\em true (resp. false) twin vertices}. In addition, a vertex $u\in V(G)$ is \emph{universal} if $N_G[u]=V(G)$. %For each $n,m\in\N$, $K_n$ and $K_{n,m}$ are, respectively, the complete graph on $n$ vertices.
%Given $U\subset V(G)$, $N(U)=\cup_{v\in U} N(v)$ and, for each $k\in \{1,\ldots, |V(G)|\}$, $\delta_k(G)=\min \{|N(U)|\,:\, U\subset V(G), |U|=k\}$. Observe that $\delta_1(G)=\delta(G)$.
If there is no ambiguity, we omit the subscript. Given a graph $G$, $\overline G$ denotes its complement. 
The subgraph induced by a set $S$ of vertices of $G$ is denoted by $G[S]$. We write $G- S$ for the subgraph induced by $V(G)- S$. If $S=\{v\}$ we simply write $G-v$.

%This section is devoted to notational and preliminary issues. For notation and graph theory terminology, %we in general follow~\cite{MHAYbookTD}. 
%A \emph{non-trivial graph} is a graph on at least two vertices. 
%For each positive integer $n$, $K_n$ and $P_n$ are, respectively, the complete graph and the path with $n$ vertices. For $n\geq 3$, $C_n$ denotes the cycle with $n$ vertices.

%Leaves are also called \emph{pendant vertices}. A \emph{strong support vertex} is a vertex with at least two leaf-neighbors. 
%We denote by $L(G)$ the set of leaves of $G$ and for each $v\in V(G)$, $L(v)=L(G)\cap N(v)$. Given $\ell \in L(G)$, $s(\ell)$ is the support vertex of $\ell$, that is, $N(\ell)=\{s(\ell)\}$.

%A vertex $v\in )$ is called a {\emph{true (resp. false) twin vertex}} if there exists $u\in V(G)-\{v\}$, such that $u$ and $v$ are true (resp. false) twin vertices.

%Given two graphs $G$ and $R$, and $v\in V(G)$, the graph obtained by \emph{replacing} $v$ \emph{by} $R$ \emph{in} $G$ is the graph whose vertex set is  $(V(G)-\{v\})\cup V(R)$ and whose edges either belong to $E(G -v)\cup E(R)$ or connect any vertex in $V(R)$ with any vertex in $N_G(v)$.

A set $D\subseteq V(G)$ is a \textit{dominating set} of $G$ if for every $v\in V(G)- D$ there exists $u \in D$ such that $uv \in E(G)$. The \textit{domination number} of $G$, denoted by $\g{G}$, is the smallest cardinality among all the dominating sets of $G$.
%(If the condition only requires that every vertex of $V(G)$ have a neighbor in $S$, then the resulting invariant is the \emph{total domination number} $\gt{G}$ of $G$.) 
Domination is one of the classical concepts in graph theory, renowned for its numerous real-world applications. Introduced over 60 years ago, it has undergone extensive development. The topic was comprehensively surveyed in recent monographs by T. Haynes, S. Hedetniemi, and M. Henning~\cite{HHH-2021,HHH-2023}. Over the years, many variations of the classical domination number of a graph have been introduced. A classical one is \emph{double domination}, which was introduced by F. Harary and T. Haynes~\cite{HH-2000} in 2000. A vertex set $D$ is a \emph{double dominating set} of $G$ if $|N[v]\cap D| \geq 2$ for every $v\in V$. The \emph{double domination number} of $G$, $\gamma_{\times2}(G)$, is the minimum cardinality of a double dominating set of $G$. Note that double domination is only defined for graphs having no isolated vertices.

Motivated by graph domination games~\cite{BHKR-2021,BKR-2010}, numerous domination-type vertex sequences have been introduced. A recent and natural one is the Grundy dominating sequence of a graph~\cite{BGMRR-2014}, which somehow describes the worst-case scenario that can occur when constructing a minimal dominating set.
%Specifically, the Grundy domination number of a graph which is the maximum length of a sequence of its vertices, such that each newly added vertex dominates at least one vertex not dominated by any preceding vertices in the sequence~\cite{BGMRR-2014}. 
Specifically, let $S=(v_1,\ldots,v_k)$ be a sequence of distinct vertices of $G$. The corresponding set $\{v_1,\ldots,v_k\}$ of vertices from the sequence $S$ will be denoted by $\widehat{S}$. The sequence $S$ is a {\em legal sequence} if
\begin{equation}
\label{e:total}
N[v_i] - \bigcup_{j=1}^{i-1}N[v_j] \ne\emptyset.
\end{equation}
holds for every $i\in\{2,\ldots,k\}$.
If, in addition, $\widehat{S}$ is a dominating set of $G$, then we call $S$ a \emph{dominating sequence} of $G$. The maximum length of a dominating sequence in $G$ is called the \emph{Grundy domination number} of $G$ and denoted by $\gr{G}$; the corresponding sequence is called a \emph{Grundy dominating sequence} of $G$. %It is easy to see that $\g{G}\leq\gr{G}$. 
The {\sc{Grundy Domination}} problem consists in finding a Grundy dominating sequence of a given graph $G$. Regarding some computational complexity results, efficient algorithms that solve this problem for trees, cographs, split graphs~\cite{BGMRR-2014}, and chain graphs~\cite{BSP-2023}. The decision version of this problem is NP-complete even when restricted to chordal graphs~\cite{BGMRR-2014}, bipartite graphs, and co-bipartite graphs~\cite{BSP-2023}.

Following this approach, T. Haynes and S. Hedetniemi~\cite{HH-2021} introduced the \emph{Grundy double dominating sequences}. If $u\in N_G[v]$, we say that $v$ \emph{dominates} $u$ and $u$ \emph{is dominated} by $v$.
%defined as a sequence of vertices such that each newly added vertex dominates at least one vertex that was dominated at most once by the preceding vertices in the sequence. 
A sequence $S =(v_1, v_2,\ldots, v_n)$ is called a \emph{double neighborhood sequence (DNS)} of $G$
if for each $i$, the vertex $v_i$ dominates at least one vertex $u$ of $G$ that is dominated
at most once by the vertices $v_1, v_2,\ldots,v_{i-1}$. If $\widehat S$ is a double dominating set of $G$, then we call $S$ a \emph{double dominating sequence (DDS)} of $G$. A double dominating sequence
of $G$ with maximum length is called a \emph{Grundy double dominating sequence (GDDS)} of
$G$. The length of a GDDS is the \emph{Grundy double
domination number} of $G$ and is denoted by $\dgr{G}$. Given a graph $G$ without
isolated vertices, the {\sc{Grundy Double Domination (GDD)}} problem asks to find
a GDDS of $G$. It is known that the decision version of this problem is NP-complete even restricted to split graphs~\cite{BPS-2021}, bipartite graphs and co-bipartite graphs~\cite{PS-2024}. On the other hand, linear-time algorithms are presented to solve the GDD problem for threshold graphs~\cite{BPS-2021} and chain graphs~\cite{PS-2024}.

In this paper, we investigate GDDSs in graphs from both a computational complexity and a structural point of view. In Section~\ref{sec:bounds}, we provide general tight bounds for the GDDN of a graph. In Section~\ref{sec:atomic}, we describe GDDSs for vertex-removed graphs, and, in particular, we obtain GDDSs of a graph $G$ in terms of GDDSs of $G-v$ , for $v$ pendant, universal, true twin, or false twin vertex. In addition, in Section~\ref{sec:modular} we provide GDDSs for the join of two graphs. Furthermore, applying these results, we prove that the GDD problem is linear for $P_4$-tidy graphs, thereby solving an open problem previously posed for cographs by B.~Bre\v{s}ar et al. in~\cite{BPS-2021}. %Finally, we characterize the graphs with $\dgr{G}\in\{2,3\}$.

% \begin{center}
% \fbox{\parbox{0.87\linewidth}{\noindent
% {\sc Grundy Domination Number Problem}\\[.8ex]
% \begin{tabular*}{.96\textwidth}{rl}
% {\em Input:} & $G=(V,E)\;,\; k\in\mathbb{Z}^+$.\\
% {\em Question:} & \text{Is there a Grundy dominating sequence of} $G$ \text{of length at least} $k$?
% \end{tabular*}
% }}
% \end{center}

% \begin{center}
% \fbox{\parbox{0.87\linewidth}{\noindent
% {\sc Grundy Double Domination Number Problem}\\[.8ex]
% \begin{tabular*}{.96\textwidth}{rl}
% {\em Input:} & $G=(V,E)\;,\; k\in\mathbb{Z}^+$.\\
% {\em Question:} & \text{Is there a Grundy double dominating sequence of} $G$ \text{of length at least} $k$?
% \end{tabular*}
% }}
% \end{center}

We conclude this section with several useful observations and definitions. 
Given $S=(v_1,\ldots,v_k)$ a sequence of distinct vertices of $G$, we denote by $o_S(v_i)$ the order of $v$ in $S$, that is, $o_S(v_i)=i$. For $R\subseteq V(G)$, $S-R$ is the subsequence obtained from $S$ by deleting the vertices in $\widehat{S}\cap R$.
Let $S_1=(v_1,\ldots , v_k)$ and $S_2=(u_1,\ldots , u_m)$, $k,m \geq 0,$ be two sequences of distinct vertices of $G$, with $\widehat S_1\cap \widehat S_2=\emptyset$. The {\em concatenation} of $S_1$ and $S_2$ is defined as the sequence $S_1 \oplus S_2=(v_1,\ldots , v_k,u_1,\ldots , u_m)$. %Clearly $\oplus$ is an associative operation on the set of all sequences, but is not commutative.
In addition, we consider the \emph{empty sequence}, denoted by $S=()$. Then we have $()\oplus S= S \oplus ()= S$, for any sequence $S$.

Given $S=(v_1,\ldots,v_k)$ a DNS of $G$, we define $N_S^1[v_i]=N_G[v_i]-\cup_{j=1}^{i-1}N_G[v_i]$, $N_S^2[v_i]=\{w\in N_G[v_i]:\sum_{j=1}^{i-1}|\{w\}\cap N_G[v_j]|=1\}$, and $N_S[v_i]=N_S^1[v_i]\cup N_S^2[v_i]$. We denote by $S^1$ and $S^2$ the subsequences of $S$ such that $\widehat{S^1}=\{v\in\widehat S:\ N_S^1[v]\neq\emptyset\}$ and $\widehat{S^2}=\widehat S- \widehat{S^1}$. Note that if $S$ is a DDS of a graph $G$ without isolated vertices, then $\widehat{S^1}\neq\emptyset$ and $\widehat{S^2}\neq\emptyset$ since $v_1\in \widehat{S^1}$ and $v_k\in \widehat{S^2}$.

%%%%%%%%%%%%%%%%%%%%%%%%%%%%%%%%%%%%%%%%%%%%%%%%%%%%%%%%%%%%%%5
\section{General bounds}\label{sec:bounds}

In this section, we consider graphs without isolated vertices. Let $S$ be a DNS of $G$. For each $u\in \widehat{S^1}$, let $P_S(u)=\{v\in\widehat{S^2}:\ N_S^1[u]\cap N_S^2[v]\neq\emptyset\}$. Note that $\cup_{u\in\widehat{S^1}} P_S(u)=\widehat{S^2}$ (some sets $P_S(u)$ may be empty) and $o_S(u)<o_S(v)$ for every $v\in P_S(u)$. Then, if $S$ is a DDS of $G$, as we have mentioned $\widehat{S^2}\neq\emptyset$ and there exists $u\in\widehat{S^1}$ such that $P_S(u)\neq\emptyset$.

Furthermore, given a sequence $S$ of distinct vertices of $G$ with $u,v\in\widehat{S}$ and $o_S(u)<o_S(v)$, we denote by $S_{u\rightarrow v}$ to the sequence obtained from $S$ by moving $u$ immediately after $v$. Formally, $\widehat{S_{u\rightarrow v}}=\widehat{S}$ and 
\[o_{S_{u\rightarrow v}}(x)=\left\{\begin{array}{ll}
    o_S(x) & \text{if }o_S(x)<o_S(u)\text{ or }o_S(x)>o_S(v), \\
    o_S(x)-1 & \text{if }o_S(u)<o_S(x)\leq o_S(v),\\
    o_S(v) & \text{if }u=v.
\end{array}
\right.\]

\begin{lemma}\label{prop:s_u}
    Let $S$ be a DDS of a graph $G$ and $u\in\widehat{S}$ such that $P_S(u)\neq\emptyset$. Consider $v_u\in P_S(u)$, the sequence $S_{u\rightarrow v_u}$, and $P'=\{v\in P_S(u):\ o_S(v)\leq o_S(v_u)\}$. Then:
        \begin{enumerate}
        \item $S^1$ is a dominating sequence of $G$.
        \item $S_{u\rightarrow v_u}$ is a DDS and $\widehat{S_{u\rightarrow v_u}}=\widehat{S}$.
        \item If $S$ is a GDDS of $G$, $S_{u\rightarrow v_u}$ also is.
        \item $P'\cup (\widehat{S^1}-\{u\})\subseteq\widehat{S_{u\rightarrow v_u}^1}$.
        \end{enumerate}
\end{lemma}
\begin{proof}
 It is easy to see that $S^1$ is a dominating sequence of $G$.
    The fact that $\widehat{S_{u\rightarrow v_u}}=\widehat{S}$ is immediate from the definition of $S_{u\rightarrow v_u}$. Let us see that $S_{u\rightarrow v_u}$ is a DDS of $G$. Observe the following:        
        \begin{itemize}
            \item If $x\in\widehat{S_{u\rightarrow v_u}}$ and $o_{S_{u\rightarrow v_u}}(x)>o_{S_{u\rightarrow v_u}}(u)$, $\emptyset\neq N_S[x]=N_{S_{u\rightarrow v_u}}[x]$.
            \item If $x\in\widehat{S_{u\rightarrow v_u}}$ and $o_{S_{u\rightarrow v_u}}(x)<o_{S_{u\rightarrow v_u}}(u)$, $\emptyset\neq N_S[x]\subseteq N_{S_{u\rightarrow v_u}}[x]$.
            \item $\emptyset\neq N_S^1[u]\cap N_S^2[v_u]\subseteq N_{S_{u\rightarrow v_u}}^2[u]$.
        \end{itemize}
        These facts show that $S_{u\rightarrow v_u}$ is a DDS of $G$.\\
        The third item follows immediately from the second.\\
        In addition, every $x\in P'$ verifies that $\emptyset\neq N_S^1[u]\cap N_S^2[x]=N_{S_{u\rightarrow v_u}}^1[x]$. Besides, $\widehat{S^1}-\{u\}\subset\widehat{S_{u\rightarrow v_u}^1}$. Therefore, $P'\cup (\widehat{S^1}-\{u\})\subseteq\widehat{S_{u\rightarrow v_u}^1}$.
\end{proof}

Let us see some bounds for $\dgr{G}$ with respect to the closely related parameters $\gamma_{\times2}(G)$ and $\gamma_{gr}(G)$.

% \begin{prop}
%     Let $x,v\in V(G)$ such that $N_G[x]\subseteq N_G[v]$. If there exists a DDS $S$ of $G$ such that $x,v\in\hat S$ and $o_S(x)>o_S(v)$. Then, $N_S^2[x]\subseteq N_S^1[v]$ and the sequence obtained by moving $v$ immediately after (or before) $x$ is a DDS of $G$.
% \end{prop}

% \begin{lemma}
%     Let $S$ be a dds of $G$, such that $x,v\in\hat S$, $x$ appears in $S$ after $v$ and $N_S[x]\subseteq N_S[v]$. Then, $N_S^2[x]\subseteq N_S^1[v]$ and the sequence obtained by moving $v$ immediately after (or before) $x$ is a dds of $G$.
% \end{lemma}

% \begin{coro}
%     Let $S$ be a dds of $G$ such that $x,v\in\hat S$ and $x$ appears in $S$ after $v$. Consider that $x$ is the last vertex in $S$ such that $N_S[x]\subseteq N_S[v]$. Then, the sequence obtained by moving $v$ immediately after (or before) $x$ is a dds of $G$.
% \end{coro}

% \begin{lemma}
%     Let $S$ be a Grundy double dominating sequence of $G$ and $u\in S^1$. Then, there exists a sequence $S'$ such that $\widehat S=\widehat S'$ and 
% \end{lemma}

\begin{prop}\label{prop:bounds}
Let $G$ be a graph without isolated vertices. Then:
\begin{enumerate}
    \item $\dfrac{2|V(G)|}{\Delta(G)+1}\leq\gamma_{\times2}(G)\leq\dgr{G}\leq |V(G)|+1-\delta(G)$,
    %\item $\gamma^t_{gr}(G)\leq\dgr{G}$,
    \item $\gr{G}+1\leq\dgr{G}\leq2\gr{G}$.
\end{enumerate}
\end{prop}
\begin{proof}
    \begin{enumerate}
        \item Notice that the sequence of vertices of a minimal double dominating set of $G$ in any order is a double dominating sequence of $G$. Then, by the known spherical bound \cite{HH-2000} we have $\frac{2|V(G)|}{\Delta(G)+1}\leq\gamma_{\times2}(G)\leq\dgr{G}$.
        The inequality $\dgr{G}\leq |V(G)|+1-\delta(G)$ follows from the fact that any set $S$ with $|V(G)|+1-\delta(G)$ vertices in $G$ verifies that $|N[v]\cap S|\geq2$.
        %\item Let $S=(v_1,\ldots,v_n)$ be a total dominating sequence of $G$. For $i=1,\dots,n$, if $u\in N(v_i) - \cup_{j=1}^{i-1}N(v_j)$, then $u\notin \cup_{j=1}^{i-1}N(v_j)$. Thus, $u\in N_S^2[v_i]$ if $u\in\{v_1,\ldots,v_{i-1}\}$ and $u\in N_S^1[v_i]$ o.w.         Therefore, $S$ is a DNS of $G$ and the result follows.
        \item To prove the lower bound, let us consider $S=(v_1,\ldots,v_n)$ a dominating sequence of $G$ and $u\in N[v_n] - \cup_{j=1}^{n-1}N[v_j]$. If $u\neq v_n$ (and so $u\notin\widehat S$), then $S\oplus (u)$ is a DNS of $G$. If $u=v_n$, since $|N[u]|\geq2$ and $N(u)\cap\widehat S=\emptyset$, there exists $x\in N(u)-\widehat S$ and $S\oplus (x)$ is a DNS of $G$. Thus, $\gr{G}+1\leq\dgr{G}$.\\
        Now, let us show that $\dgr{G}\leq2\gr{G}$. Let $S$ be a GDDS of $G$ that maximizes $|\widehat{S^1}|-|\widehat{S^2}|$. From Lemma~\ref{prop:s_u}, the subsequence $S^1$ is a dominating sequence of $G$.\\
        Let us assume that $|\widehat{S^1}|- |\widehat{S^2}|< 0$, that is, $|\widehat{S^1}|<|\widehat{S^2}|$. Since $|\widehat{S^1}|<|\widehat{S^2}|$, there exists $u\in\widehat{S^1}$ such that $|P_S(u)|\geq2$. Let $v_u\in P_S(u)$ with maximum order in $S$. Consider $S_{u\rightarrow v_u}$ and $P'$ as in Lemma~\ref{prop:s_u}. Observe that $P'=P_S(u)$. Again from Lemma~\ref{prop:s_u}, $S_{u\rightarrow v_u}$ is a GDDS of $G$ and $P_S(u)\cup (\widehat{S^1}-\{u\})\subseteq\widehat{S_{u\rightarrow v_u}^1}$. Since $P_S(u)\cap (\widehat{S^1}-\{u\})=\emptyset$ and $|P_S(u)|\geq2$, then $|\widehat{S_{u\rightarrow v_u}^1}|>|\widehat{S^1}|$ and $|\widehat{S_{u\rightarrow v_u}^2}|<|\widehat{S^2}|$. Therefore, $S_{u\rightarrow v_u}$ is a GDDS and $|\widehat{S_{u\rightarrow v_u}^1}|- |\widehat{S_{u\rightarrow v_u}^2}|>|\widehat{S^1}|- |\widehat{S^2}|$, which contradicts the election of $S$.\\
        Thus, $|\widehat{S^1}|- |\widehat{S^2}|\geq 0$ and then $\dgr{G}=|\widehat{S^1}|+|\widehat{S^2}|\leq2|\widehat{S^1}|\leq2\gr{G}$.\\
    \end{enumerate}
\end{proof}

Notice that for complete graphs with at least two vertices all the bounds presented above are tight. In addition, from the proof of item $2$ of Proposition~\ref{prop:bounds} we have the following.

\begin{coro}\label{coro:grundyS1}
    There exists a GDDS $S$ of $G$ such that $|\widehat{S^1}|\geq |\widehat{S^2}|$. Furthermore, $|\widehat{S^1}|\geq\frac{\dgr{G}}{2}$.
\end{coro}

Given a graph parameter, an interesting problem is characterizing the graphs for which the mentioned parameter is small. Let $|V(G)|=n$. It is not hard to see that $\dgr{G}=2$ if and only if $G$ is isomorphic to $ K_n$ ($n\geq2$). Assume that $\dgr{G}=3$ and $deg_G(v)\leq n-3$. Let $v_1,v_2\notin N_G[v]$ and $v_3\in N_G(v)$. It is easy to see that $(v_1,v_2,v_3,v)$ is a DNS. Thus, $deg_G(v)\geq n-2$ for all $v\in V(G)$. We also know that $G$ is not a complete graph. Therefore, $G$ is isomorphic to $\oplus_{i=1}^kG_i$ with $k\geq2$, $G_1=\overline{K_2}$, and for every $i=2,\dots,k$, $G_i\in\{\overline{K_2},K_1\}$.
\begin{rem}
    Let $G$ be a connected graph without isolated vertices. Then,
    \begin{itemize}
        \item If $\dgr{G}=2$, $G$ is a complete graph with at least two vertices.
        \item If $\dgr{G}=3$,  $G$ is isomorphic to $\vee_{i=1}^kG_i$ with $k\geq2$, $G_1=\overline{K_2}$, and for every $i=2,\dots,k$, $G_i\in\{\overline{K_2},K_1\}$.
    \end{itemize}
\end{rem}

% Regarding extreme values of $\dgr{G}$, from the upper bound in item $1$ of Proposition~\ref{prop:bounds}, if $\dgr{G}=n$, then $G$ has no pendant vertices.
% In addition, by~\ref{prop:pendant}, if every vertex of $G$ is either a pendant vertex or adjacent to a pendant vertex, then $\dgr{G}=n$. In this sense, an interesting problem is characterizing connected graphs with $\dgr{G}=n$. Examples of graphs with $\dgr{G}=n$ seen in this paper are trees and thin spider graphs (S,C,$\emptyset$). 

%%%%%%%%%%%%%%%%%%%%%%%%%%%%%%%%%%%%%%%%%%%%%%%%%%%%%%%%%%%%%%%%
\section{GDDSs and vertex deletion}\label{sec:atomic}

Studying how graph parameters behave when vertices are deleted/added is a useful tool, for instance, for developing dynamic algorithms to determine their values. To this end, we introduce a slightly generalization of the concept of GDDS of graphs.

Given a graph $G$, a sequence $S$ is a \emph{maximum double neighborhood sequence} (MDNS) of $G$ if $S$ is a DNS of $G$ with maximum length. The length of a MDNS is the \emph{maximum double neighborhood number} of $G$ and we denote it by $\dgri{G}$. Recall that $I(G)$ is the set of isolated vertices of $G$. If $G$ has at least one edge and a sequence $S$ of $G$ is a MDNS of $G$ then $I(G)\subset\widehat{S}$ and $S-I(G)$ is a GDDS of $G- I(G)$. In addition, $\dgri{G}=\dgr{G- I(G)}+|I(G)|$ if $E(G)\neq\emptyset$ and $\dgri{G}=|I(G)|$, otherwise. Notice that if $G$ has no isolated vertices, then $S$ is a MDNS of $G$ if and only if $S$ is a GDDS of $G$, and so $\dgri{G}=\dgr{G}$. Consider the graph function $a$ such that $a(G)=1$ if $I(G)\neq\emptyset$ and $a(G)=0$, otherwise.

Our next result shows that $\dgr{G}$ does not increase upon vertex deletion and can decrease by a maximum of $3$ units.

%\begin{rem}
%    Let $G$ be a graph without isolated vertices, $S$ a GDDS of $G$, and $v\in V(G)$. Then, $|N_G[v]\cap\hat S|\geq2$.
%\end{rem}

\begin{prop}\label{prop:dgrG-v}
    Let $v\in V(G)$ such that $G-v$ has no isolated vertices. Then, 
    \[\dgri{G-v}\leq\dgri{G}\leq\dgri{G-v}+3.\] 
\end{prop}
\begin{proof}
    Notice that every DNS of $G-v$ is a DNS of $G$ and thus $\dgri{G-v}\leq\dgri{G}$.\\
    Let $S$ be a MDNS of $G$. Let $v_1,v_2\in\widehat{S}$ such that $v\in N_S^1[v_1]$ and $v\in N_S^2[v_2]$. Observe that $v_1$ or $v_2$ may be $v$. It follows that $N_S[x]-\{v\}\neq\emptyset$ for every $x\in\widehat{S}-\{v_1,v_2,v\}$. Then, the sequence $S-\{v_1,v_2,v\}$ is a DNS of $G-v$. Therefore, $\dgr{G}\leq\dgr{G-v}+3$.
\end{proof}

Observe that the bounds presented in Proposition~\ref{prop:dgrG-v} are tight. As usual, for $n,m\in\N$, let denote by $K_{n,m}$ the complete bipartite graph with bipartition $(X,Y)$ such that $|X|=n$ and $|Y|=m$. It is known that if $n>m$, $\dgr{K_{n,m+1}}=\dgr{K_{n,m}}=n+1$ \cite{PS-2024}. Now, let $G$ and $H$ the graphs in Figure~\ref{fig:propdgrG-v}. The sequence $S=(v_1,\dots,v_8,v,v_9,v_{10})$ is a GDDS of $G$ and $S'=(v_1,\dots,v_6,v_9,v_{10})$ is a GDDS of $G-v$ (connected). The sequence $S=(v_1,v_2,v_3,v_4,v_7,v_8,v,v_5,v_6)$ is a GDDS of $H$ and $S'=(v_1,\dots,v_6)$ is a MDNS of $G-v$ (not connected).
%Upper bound $G-v=K_{n,3}+K_1$ ($n\geq3$), $v$ adjacent to $K_1$ and two vertices of $S_3$.\\

\begin{figure}[h]
    \centering
    \begin{tikzpicture}
        \foreach \x in {1,2,...,4}{
        \foreach \y in {5,6,...,8}{
        \draw (-3,\x) -- (-1,\y-4);
        }}
        \foreach \x in {1,2,...,4}{
        \filldraw (-3,\x) circle (2.5pt) node[left] {$v_{\x}$};
        }
        \foreach \y in {5,6,...,8}{
        \filldraw (-1,9-\y) circle (2.5pt) node[above] {$v_{\y}$};
        }
        \foreach \x in {1,3}{
        \foreach \y in {1,3}{
        \draw (1,3) -- (3,\y);
        \filldraw (1,\x) circle (2.5pt);
        \filldraw (3,\y) circle (2.5pt);
        }}
        \draw (1,1) -- (3,1);
        \draw (1,1) -- (-1,1);
        \draw (1,1) -- (-1,2);
        \draw (1,3) -- (-1,3);
        \draw (3,3) -- (3,1);
        \node[above] at (1,1) {$v$};
        \node[above] at (3,3) {$v_{9}$};
        \node[below] at (3,1) {$v_{10}$};
        \node at (0,0) {$G$};
        %%%%%%%%%%%%%%%%%%%%%%%%%%%%%%%%%%%%%%%%%%
        \foreach \y in {1,2,...,6}{
        \filldraw (\y+4,1) circle (2.5pt) node[below] {$v_{\y}$};
        }
        \foreach \y in {7,9,11}{
        \draw (\y-1.5,2) -- (\y-2,1) -- (\y-1,1) -- (\y-1.5,2);
        \draw (\y-1.5,2) -- (7.5,3.5);
        }
        \filldraw (5.5,2) circle (2.5pt) node[left] {$v_{7}$};
        \filldraw (7.5,2) circle (2.5pt) node[left] {$v_{8}$};
        \filldraw (9.5,2) circle (2.5pt) node[left] {$v_{9}$};
        \filldraw (7.5,3.5) circle (2.5pt) node[above] {$v$};
        \node at (7.5,0) {$H$};
   \end{tikzpicture}
   \caption{$\dgr{G}=\dgr{G-v}+3$ and $\dgr{H}=\dgri{H}=\dgri{H-v}+3$.}
\label{fig:propdgrG-v}
\end{figure}

In addition, observe that for the path graph $P_n$ ($n\geq2$), $\dgr{P_{n+1}}=n+1=\dgr{P_{n}}+1=\dgr{P_{n+1}-v}+1$, where $v$ is a pendant vertex of $P_n$. Similarly, if $P$ is the paw graph, $\dgr{P}=4=\dgr{K_{3}}+2=\dgr{P-v}+2$, where $v$ is the pendant vertex of $P$.

If a vertex $v$ with particular properties is removed from a graph $G$, the result of Proposition~\ref{prop:dgrG-v} may be more specific. In what follows, we investigate the cases in which $v$ is a pendant vertex, a universal vertex, a true twin vertex, or a false twin vertex.

\begin{prop}\label{prop:pendant}
    Let $p$ be a pendant vertex of $G$, with $N(p)=\{v\}$. Then:
    \begin{enumerate}
        \item $p$ and $v$ belong to every MDNS of $G$. In addition, there exists a MDNS $S$ of $G$ such that $p$ and $v$ are the last vertices of $S$.
        \item If there exists a MDNS $S'$ of $G-p$ such that $v\notin\widehat{S'}$, then $\dgri{G}=\dgri{G-p}+2$. Else, $\dgri{G}=\dgri{G-p}+1$.
    \end{enumerate}
\end{prop}
\begin{proof}
\begin{enumerate}
    \item The result follows from the fact that if $S$ is a DNS of $G$, then $S'=(S-\{v,p\})\oplus (v,p)$ is a DNS of $G$, since $p\in N_{S'}^1[v]$ and $p\in N_{S'}^2[p]$.
    \item It is not hard to see that $\dgri{G}\leq\dgri{G-p}+2$. Note that if $S$ is a MDNS of $G$ then $S-\{v,p\}$ is a DNS of $G-p$. On the other hand, if $S'$ is a DNS of $G-v$, then $S=S'\oplus (p)$ is a DNS of $G$. Thus, $\dgri{G-p}+1\leq\dgri{G}$.\\
    If $S'$ is a MDNS of $G-p$ such that $v\notin\widehat{S'}$, then $S'\oplus (v,p)$ is a MDNS of $G$, that is, $\dgri{G}=\dgri{G-p}+2$. In addition, if $\dgri{G}=\dgri{G-p}+2$ and $S$ is a MDNS of $G$, then $S-\{v,p\}$ is a MDNS of $G-p$ such that $v\notin\widehat{S'}$.
\end{enumerate}
\end{proof}

From Proposition \ref{prop:pendant} follows that if $T$ is a tree, then $\dgri{T}=|V(T)|$, which is proved in \cite{HH-2021}.
In addition, it is easy to see that a GDDS of $T$ could be obtained by rooted $T$ at an arbitrary vertex and ordering the vertices in the sequence according to the level in increased order.

Consider now a graph $G$ such that $u$ is a universal vertex of $G$. Let $S'$ be a MDNS of $G-u$.
Assume first that $G-u$ has no isolated vertices. We have $|\widehat{S}|\geq2$ and thus $S'$ is also a DDS of $G$. Consequently, $\dgr{G}\geq\dgr{G-u}$. In addition, if $S=(v_1,\dots,v_k)$ is a DDS of $G$, then $N_S[v_i]-\{u\}\neq\emptyset$ for all $i\in\{3,\dots,4\}$ since $u\in N^1_S[v_1]$ and $u\in N^2_S[v_2]$. Furthermore, from the fact that $G-u$ has no isolated vertices follows that $N_S[v_i]-\{u\}\neq\emptyset$ for $i\in\{1,2\}$. Hence, if $u\notin\widehat{S}$, $S$ is a DDS of $G-u$. Assume that $u\in\widehat{S}$. It is not hard to see that if $o_S(u)<k$, then $v_k\in P_S(u)$ and from Lemma~\ref{prop:s_u} with $v_u=v_k$, we have that $S_{u\rightarrow v_k}$ is a DDS of $G$ such that $o_S(u)=k$. Thus, we can assume that $u$ is the last vertex of the sequence $S$. Let $x\in N^2_S[u]$. Observe that $|N_{G-u}[x]\cap\widehat{S}|=1$ and $|N_{G-u}[x]|\geq2$ ($G-u$ has no isolated vertices). Consequently, there exists a vertex $y\in N_G[x]-\widehat{S}$. Therefore, the sequence $S'=(v_1,\dots,v_{k-1},y)$ is a DDS without the vertex $u$. As we have mentioned, $S'$ is also a DDS of $G-u$. Therefore, $\dgr{G}=\dgr{G-u}$ and every MDNS of $G-u$ is a MDNS of $G$.

Finally, if $I(G-u)\neq\emptyset$, then $S=S'\oplus (u)$ is a MDNS of $G$ since $N_S[u]=N_S^2[u]=I(G)\neq\emptyset$. Then, $\dgri{G}\geq\dgri{G-u}+1$ (*).% Let $S$ be a MDNS of $G$.
%Observe that every vertex $p\in I(G-u)$ is a pendant vertex in $G$ with $N_G(p)=\{u\}$. From Proposition~\ref{prop:pendant}, $I(G-u)\cup\{u\}\subseteq\widehat{S}$. 
Let $I(G-u)=\{s_1,\dots,s_k\}$ ($k\geq1$). If $G-u$ is an edgeless graph, $G$ is isomorphic to the star graph $K_{1,|V(G)|-1}$, $S=(s_1,\dots,s_k,u)$ is a MDNS of $G$, and $\dgri{G}=|V(G)|=|I(G-u)|+1=\dgri{G-u}+1$. Otherwise, $u$ is a universal vertex in $G'=G[V(G)-I(G-u)]$ and $G'-u$ has no isolated vertices. We know that $\dgri{G-u}=\dgri{G'-u}+|I(G-u)|$. In addition, if $S$ is a MDNS of $G'-u$, then $S\oplus(s_1,\dots,s_k,u)$ is a MDNS of $G$. Thus, $\dgri{G}\leq\dgri{G'-u}+|I(G-u)|+1$. Therefore, $\dgri{G}\leq\dgri{G'-u}+|I(G-u)|+1=\dgri{G-u}+1$. The equality follows from (*).

In sum, we conclude the next result.

\begin{prop}\label{prop:universal}
    Let $G$ be a graph such that $u$ is a universal vertex of $G$. Then, 
    \[\dgr{G}=\dgr{G-u}+a(G-u).\]
    Moreover, if $a(G-u)=0$ every GDDS of $G-u$ is a GDDS of $G$. If $a(G-u)=1$ and $S$ is a MDNS of $G-u$, then $S\oplus (u)$ is a MDNS of $G$.
\end{prop}

A threshold graph $G$ is known to be constructible from a single vertex through a sequence of operations, each involving the addition of either an isolated or a universal vertex. Obtaining such a sequence, as well as the recognition of threshold graphs, can be done in linear time~\cite{threshold}. Note that if $v$ an isolated vertex of a graph $G$, then $\dgri{G}=\dgri{G-v}+1$ and if $S$ is a MDNS of $G-v$ then $S\oplus (v)$ is a MDNS of $G$. Therefore, these facts and the result in Proposition~\ref{prop:universal} led to a similar linear-time algorithm to find a MDNS of $G$ to that obtained in~\cite{BSP-2023}.
% If $G$ is connected, this process invariably concludes with the addition of a universal vertex. However, adding $k$ isolated vertices followed by a universal vertex $u$ is equivalent to first adding $u$ as a universal vertex and then attaching $k$ pendant vertices to $u$. It is also worth noting that u remains a universal vertex even after these $k$ isolated vertices are added. Therefore, $G$ can be generated by successively applying operations that either add a universal vertex or a pendant vertex adjacent to one. Finally, the result in propositions~\ref{prop:pendant} and 
\medskip

Now, let us analyze DNSs in graphs with a pair of true twin vertices.

\begin{lemma}\label{lemma:true}
    Let $v$ and $v'$ be true twin vertices of $G$ and $S$ a DNS of $G$. Then:
    \begin{enumerate}
        \item If $v'\in N_S^i[x]$, then $v\in N_S^i[x]$, for every $x\in\widehat{S}$ and $i=1,2$.
        \item If $v$ and $v'$ belong to $S$, then there exists a DNS $S'$ of $G$ such that $v$ and $v'$ appear consecutively in $S'$ (that is, $|o_S(v)-o_S(v')|=1$) and $\widehat{S'}=\widehat{S}$.
    \end{enumerate}
\end{lemma}
\begin{proof}
\begin{enumerate}
    \item It is straightforward from the definition of true twin vertices.
    \item We have $v,v'\in\widehat{S}$. Without loss of generality, we assume that $o_S(v)<o_S(v')$. Thus, $N_S^1[v]\neq\emptyset$, $N_S^1[v']=\emptyset$, and $\emptyset\neq N_S^2[v']\subseteq N_S^1[v]$. It follows that $v'\in P_S(v)$. Therefore, considering $u=v$ and $v_u=v'$ in Lemma~\ref{prop:s_u}, then $S'=S_{v\rightarrow v'}$ is a DNS of $G$ where $v$ and $v'$ appear consecutively and $\widehat{S'}=\widehat{S}$. 
\end{enumerate}
\end{proof}
%If $\dgr{G-v'}+1=\dgr{G}$ then $v$ and $v'$ belong to every GDDS of $G$.
%If there exists a GDDS $S$ of $G$ such that $v'\notin\widehat S$, then $\dgr{G-v'}=\dgr{G}$.    If $G$ is the paw graph, then $\dgr{G-v'}+1=\dgr{G}$.\\

\begin{prop}\label{prop:true}
    Let $v$ and $v'$ be true twin vertices of $G$. Then, 
    \[\dgri{G-v'}\leq\dgri{G}\leq\dgri{G-v'}+1.\]
    % if there exists a MDNS $S$ of $G$ such that $v'\notin\widehat S$, then $\dgri{G-v'}=\dgri{G}$. Moreover, 
    In addition, $\dgri{G}=\dgri{G-v'}+1$ if and only if there exists a MDNS $S'$ of $G'$ such that $v\in\widehat{S'}$, $N_{S'}^1[v]\neq\emptyset$, and for every $x\in\widehat S-\{v\}$, $N_{S'}[x]-N_{S'}^1[v]\neq\emptyset$.
\end{prop}
\begin{proof}
    We denote $G'=G-v'$. Let us see that $\dgri{G}\leq\dgri{G'}+1$. Let $S$ be a MDNS of $G$ and $S'=S-\{v'\}$. Observe that for every $x\in\widehat{S'}$, $N_S[x]=N_{S'}[x]-\{v'\}$. Then, from Lemma~\ref{lemma:true}, for every $x\in\widehat{S'}$, $N_{S'}[x]\neq\emptyset$ and $S'$ is a DNS of $G'$ with $\dgri{G'}\geq|\widehat{S'}|\geq\dgri{G}-1$. The fact that $\dgri{G'}\leq\dgri{G}$ follows from Proposition~\ref{prop:dgrG-v}.
    
    Let us see now that $\dgri{G}=\dgri{G'}+1$ if and only if there exists a MDNS $S'$ of $G'$ such that $v\in\widehat{S'}$, $N_{S'}^1[v]\neq\emptyset$, and $N_{S'}[x]-N_{S'}^1[v]\neq\emptyset, \forall x\in\hat S-\{v\}$.

    Let $S'$ be a MDNS of $G'$ such that $v\in\widehat{S'}$, $N_{S'}^1[v]\neq\emptyset$, and $N_{S'}[x]-N_{S'}^1[v]\neq\emptyset, \forall x\in\hat S-\{v\}$. Consider the sequence $S$ in $G$ obtained from $S'$ by adding $v'$ immediately after $v$. Formally, $S$ is the sequence in $G$ such that $\widehat{S}=\widehat{S'}\cup\{v'\}$, $o_S(v')=o_{S'}(v)+1$; if $o_{S'}(x)\leq o_{S'}(v)$, $o_{S}(x)=o_{S'}(x)$; and if $o_{S'}(x)>o_{S'}(v)$, $o_{S}(x)=o_{S'}(x)+1$. Then:
    \begin{itemize}
        \item $N_S[v']=N_S^2[v']=N_{S'}^1[v]\neq\emptyset$.
        \item $\emptyset\neq N_S[x]=N_{S'}[x]\neq\emptyset$ for every $x\in\widehat S$ with $o_S(x)\leq o_S(v)$.
        \item $N_S[x]=N_{S'}[x]-N_S[v']=N_{S'}[x]-N_{S'}^1[v]\neq\emptyset$ for every $x\in\widehat S$ with $o_S(x)>o_S(v')=o_S(v)+1$.
    \end{itemize}
    Thus, $S$ is a MDNS of $G$ and $|\widehat{S}|=|\widehat{S'}|+1$. From the previous bound follows that $\dgri{G}=\dgri{G'}+1$.
    
    Conversely, assume now that $\dgri{G}=\dgri{G'}+1$. Applying similar arguments as in the first part of this proof, if $S$ is a MDNS of $G$ such that $v'\notin\widehat S$, then $S'=S-v'=S$ is a DNS of $G'$ and $\dgri{G'}\geq|\widehat{S'}|\geq\dgri{G}$. Thus, $\dgri{G'}=\dgri{G}$. Then, we have that $v$ and $v'$ belong to every MDNS of $G$. From Lemma~\ref{lemma:true}, there exists a MDNS $S$ of $G$ such that $o_S(v')=o_S(v)+1$. Consider $S'=S-\{v'\}$. Thus, $S'$ is a MDNS of $G'$, $v\in\widehat{S'}$, $N_{S'}^1[v]\neq\emptyset$, and $\emptyset\neq N_{S}[x]-(N_{S}^1[v]\cup N_{S}^2[v'])\subseteq N_{S'}[x]-N_{S'}^1[v], \forall x\in\widehat S-\{v\}$.
\end{proof}

Given a graph $G$ and a function $f:V(G)\mapsto\N$, we denote $G_f$ the graphs obtained from $G$ by adding exactly $f(v)$ true twin vertices to each $v\in V(G)$, that is, $V(G_f)=\cup_{v\in V(G)}\{v^1,\dots,v^{f(v)+1}\}$ and $v^iu^j\in E(G_f)$ if and only if $vu\in E(G)$. For each $v\in V(G)$, let $G_f(v)=\{v^1,\dots,v^{f(v)+1}\}$. Note that if $f(v)=k$ for all $v\in V(G)$, then $G_f$ is isomorphic to the lexicographic product $G\circ K_{k+1}$.

% \begin{prop}
%     Let $v,v'\in V(G)$ with $v$ and $v'$ true twins in $G$. Consider $S$ a dds of $G$ such that $v,v'\in\hat S$. There exists a dds $S'$ of $G$ such that $v$ and $v'$ appear consecutively and $\hat S=\hat{S'}$.
% \end{prop}
    %If $W$ is a set of true twins in $G$ and $S$ is a DDS of $G$, then $|W\cap\widehat S|\leq2$. In addition, 
Let us see that $2\gr{G}\leq\dgr{G_f}$. Consider $S=(v_1,\dots,v_k)$ a Grundy dominating sequence of $G$ and the sequence $S'=(v^1_1,v^2_1,v_2^1,v^2_2,\dots,v^1_k,v^2_k)$ in $G_f$. Observe that for $i=2,\dots,k$, $N_{S'}^1[v^1_i]=\cup_{v\in N[v_i]-\cup_{j=1}^{i-1}N[v_j]}G_f(v)\neq\emptyset$, since $S$ is a Grundy dominating sequence of $G$. Besides, for $i=1\dots,k$, $N_{S'}[v_i^2]=N_{S'}^2[v_i^2]=N_{S'}^1[v_i]\neq\emptyset$. Therefore, $S'$ is a DDS of $G_f$ and the next remark follows.

\begin{rem}\label{rem:truetwin1}
For every graph $G$, $2\gr{G}\leq\dgr{G_f}$.
\end{rem}

\begin{lemma}\label{lemma:GfG}
    Let $G$ be a graph and a function $f:V(G)\mapsto\N$. Then, $\dgr{G_f}\leq2\gr{G}$. Moreover, there exists a GDDS $S=(v_1,\dots,v_{2k})$ of $G_f$ such that for $i=1,3,\dots,2k-1$, $v_i$ and $v_{i+1}$ are true twin vertices in $G$. 
\end{lemma}
\begin{proof}
    From Corollary~\ref{coro:grundyS1}, let $T$ be a GDDS of $G_f$ such that $|\widehat{T^1}|\geq |\widehat{T^2}|$. By Lemma~\ref{prop:s_u}, $T^1$ is a dominating sequence of $G_f$. Observe that there are no pairs of true twin vertices in $T^1$. Without loss of generality, we assume that $T^1=(v^1_1,v_2^1,\dots,u^1_k)$. It is not hard to see that $(v_1,v_2,\ldots,v_k)$ is a dominating sequence of $G$. Then, $\dgr{G_f}\leq 2|\widehat{T^1}|\leq2\gr{G}$.
\end{proof}

As a by-product of Remark~\ref{rem:truetwin1} and Lemma~\ref{lemma:GfG} we have the following.

\begin{prop}\label{2G=Gf}
    Let $G$ be a graph and a function $f:V(G)\mapsto\N$. Then, $2\gr{G}=\dgr{G_f}$.
\end{prop}

We said that a graph class $\mathcal{F}$ closed under addition of true twin vertices if for every $G\in\mathcal F$ and every $v\in V(G)$, the graph obtained from $G$ by adding a true twin vertex to $v$ also belongs to $\mathcal{F}$. Notice that Proposition~\ref{2G=Gf} provide a reduction of {\sc Grundy Domination} for $\mathcal{F}$ into {\sc Grundy Double Domination} for $\mathcal{F}$.

\begin{coro}\label{coro:complexity}
 Let $\mathcal{F}$ be a graph class closed under the addition of true twin vertices. If the decision version of {\sc Grundy Domination} is NP-complete for $\mathcal F$, then the decision version of {\sc Grundy Double Domination} also is NP-complete for $\mathcal{F}$.    
\end{coro}

As we have mentioned, the decision version of {\sc Grundy Double Domination} is NP-complete for chordal graphs~\cite{BSP-2023} and for co-bipartite graphs~\cite{PS-2024}. The above corollary provides an alternative proof of these results since the decision version of {\sc Grundy Domination} is NP-complete for chordal graphs~\cite{BGMRR-2014} and for co-bipartite graphs~\cite{BPS-2021}.

%Calcular grundy en cordales es NP-c. Como $G_f$ es cordal para todo $G$ cordal, tenemos que ddg en cordales es NP-c.

We end this section providing a bounds for the length of a MDNS of a graph when a false twin vertex is removed. 

\begin{prop}\label{prop:false}
    Let $v$ and $v'$ false twin vertices of $G$. Then,
    \[\dgri{G-v'}\leq\dgri{G}\leq\dgri{G-v'}+1.\] 
    In addition, if $S$ is a MDNS of $G$ and $\dgri{G}=\dgri{G-v'}+1$, then $\widehat{S}\cap\{v,v\}\neq\emptyset$.
\end{prop}
\begin{proof}
    If $N_G(v)=N_G(v')=\emptyset$, both vertices are isolated and the result is trivial. Assume $N_G(v)=N_G(v')\neq\emptyset$.
    Consider $G'=G-v'$. From Proposition~\ref{prop:dgrG-v}, we have $\dgr{G-v'}\leq\dgr{G}$. %Let $G'=G-v'$ and $S'$ be a double neighborhood sequence of $G-v'$. If $|N_{G'}[v]\cap\widehat{S'}|\geq2$ then $S'$ is also a double neighborhood sequence of $G$, else $S'\oplus (v')$ is a double neighborhood sequence of $G$. Thus, $\dgr{G-v'}\leq\dgr{G}$.

    In order to prove the remaining inequality, consider $S$ a MDNS of $G$. If $v,v'\notin\widehat{S}$, then for every $x\in\widehat{S}$, $v'\in N_S[x]$ if and only if $v\in N_S[x]$. Hence, $S'=S$ is also a MDNS of $G'$.

    Assume now that $v\in\widehat{S}$ and $v'\notin\widehat{S}$ (analogously if $v'\in\widehat{S}$ and $v\notin\widehat{S}$). Then, there exist $v_1,v_2\in\widehat{S}\cap N_G(v)$ with $v_1\neq v_2$ such that $v'\in N_S^1[v_1]$ and $v'\in N_S^2[v_2]$. Note that $v_1$ and $v_2$ are the two lowest-order vertices of $S$ in $N_G(v)$. If $o_S(v_2)<o_S(v)$, then $v\in N_S^1[v_1]$ and $v\in N_S^2[v_2]$. Hence, $S$ is also a MDNS of $G'$. If $o_S(v)<o_S(v_2)$, we have $v\in N_S[v_1]$. Thus, $S'=S-\{v_2\}$ is a DNS of $G'$.

    Finally, assume that $v,v'\in\widehat{S}$ and $o_S(v)<o_S(v')$ (analogously if $o_S(v')<o_S(v)$).
    
    If $v'\in N_S^1[v']$, then $v\in N_S^1[v]$ and $v,v'\in N_S^2[u]$ with $u$ the lowest-order vertex of $S$ in $N_G(v)$. Hence, $N_S[u]-\{v'\}\neq\emptyset$ and for every $x\in\widehat{S}-\{u,v'\}$, $v'\notin N_S[x]$. Thus, $S'=S-\{v'\}$ is a DNS of $G'$.

    Consider $v'\notin N_S^1[v']$ and let $v_1,v_2\in\widehat{S}$ with $v'\in N_S^1[v_1]$ and $v'\in N_S^2[v_2]$. Observe that $v_1\in N_G(v)$ and $v\in N_S[v_1]$. If $v_2=v'$ or $N_S[v_2]\neq\{v'\}$, then $S'=S-\{v'\}$ is a DNS of $G'$. Otherwise, $v_2\neq v'$ and $N_S[v_2]=\{v'\}$. Let $v_1\in\widehat{S}\cap N_G(v)$ with $v_1\neq v_2$ be such that $v'\in N_S^1[v_1]$. Hence, $o_S(v_1)<o_S(v_2)<o_S(v')$, $N_S[v']=N_S^2[v']$ and $v',v_1,v_2\notin N_S[v']$. Let $u\in N_S[v']$. Note that $u\in N_G(v)$. Consider the sequence $S'$ obtained from $S$ by moving $v$ to the position of $v'$ and deleting $v'$, that is, $S'=S_{v\rightarrow v'}-\{v'\}$. We claim that $S'$ is a DNS of $G'$. Note that $v\in N_{S'}[v_1]$, $v\in N_{S'}[v_2]$, $u\in N_{S'}[v]$ and for every $x\in\widehat{S'}-\{v,v_1,v_2\}$, $\emptyset\neq N_S[x]\subseteq N_{S'}[x]$, since $v_2\notin N_S[x]$. Therefore, $S'$ is a DNS of $G'$.

    In sum, for any case, the sequence $S'$ is a DNS of $G'$ and $|\widehat{S'}|\geq |\widehat{S}|-1=\dgri{G}-1$. Thus, $\dgri{G}\leq\dgri{G'}+1$.
    
    In addition, as we have mentioned, if $v,v'\notin\widehat{S}$ then $S'=S$ is also a MDNS of $G'$. In consequence, $\dgri{G}=\dgri{G-v'}+1$, then $\widehat{S}\cap\{v,v\}\neq\emptyset$.
\end{proof}

% PENSAR EN COND. SUF Y/O NEC. VER QUE PASA EN BIP. DIST. HERED.

% \begin{prop}
%     Let $v$ be a simplicial vertex of $G$. Then, there exists a maximum dds $S$ of $G$ such that $v\in\hat S$ and $v$ appears in $S$ after the first of its neighbors.
% \end{prop}

% \begin{prop}
%     Let $v,v'$ two simplicial false twin vertices of $G$. Then, $\dgr{G}=\dgr{G-v'}+1$. Besides, there exists a maximum dds of $G$ such that $v,v'\in\hat S$ and $v,v'$ appear consecutively in $S$.  
% \end{prop}

%%%%%%%%%%%%%%%%%%%%%%%%%%%%%%%%%%%%%%%%%%%%%%
\section{MDNSs under modular decomposition}\label{sec:modular}

A common strategy for studying various problems in graph theory is to decompose a graph into simpler, smaller components. While numerous decomposition methods exist and have been applied to a wide array of problems, this work utilizes the modular decomposition of graphs to study DNSs. This approach has proven valuable for analyzing problems in graph families characterized by their $P_4$ structure, including cographs, $P_4$-tidy graphs, partner limited graphs, interval graphs, and permutation graphs; see \cite{CH-1994} and the references therein.

To this end, let us consider the graph operations \emph{disjoint union} and \emph{join}. Given non-empty graphs $G$ and $H$ with disjoint sets of vertices, $G+H$ denotes their disjoint union and $G\vee H$ their join. Let us first study Grundy double dominating sequences for the disjoint union and join of graphs. Observe also that $\dgri{G\vee H}=\dgr{G\vee H}$.

\begin{lemma}\label{lemma:unionyjoin}
    Let $G$ and $H$ be graphs. Consider $S_G$ and $S_H$ MDNS of $G$ and $H$, respectively. Then,
    \begin{enumerate}
        \item $\dgri{G+H}=\dgri{G}+\dgri{H}$ and $S_G\oplus S_H$ is a MDNS of $G+H$.
        \item $\dgri{G\vee H}=\max\{\dgri{G}+a(G),\dgri{H}+a(H)\}$. In addition:
        \begin{enumerate}
            \item If $\dgri{G}+a(G)\geq\dgri{H}+a(H)$ and $a(G)=0$, then $S_G$ is a MDNS of $G\vee H$.
            \item If $\dgri{G}+a(G)\geq\dgri{H}+a(H)$ and $a(G)=1$, then, for any $h\in V(H)$, $S_G\oplus (h)$ is a MDNS of $G\vee H$.
            \item If $\dgri{G}+a(G)<\dgri{H}+a(H)$ and $a(H)=0$, then $S_H$ is a MDNS of $G\vee H$.
            \item If $\dgri{G}+a(G)<\dgri{H}+a(H)$ and $a(G)=1$, then, for any $g\in V(G)$, $S_H\oplus (g)$ is a MDNS of $G\vee H$.
        \end{enumerate}
    \end{enumerate}
\end{lemma}
\begin{proof}
    It is straightforward that $S_G\oplus S_H$ is a MDNS of $G+H$ and $\dgri{G+H}=\dgri{G}+\dgri{H}$.
    %Now, let us consider the graph $G\oplus H$ and assume that $|\widehat{S_G}|+a(G)\geq |\widehat{S_H}|+a(H)$.
    %Note that if $|V(H)|=1$ with $V(H)=\{h\}$, then the results follows from the fact that a MDNS of $G\oplus H$ is $S_G$ if $a(H)=0$ and $S_G\oplus (h)$ if $a(G)=1$. If $|V(G)|=1$, then $|V(H)|=1$ and the result follows.\\
    %Assume that $G$ and $H$ have both at least two vertices.\\
    %Let us first see that $\dgri{G\oplus H}\geq\dgri{G}+a(G)$. Let $h\in V(H)$. It is not hard to see that $S_G$ is a double dominating sequence of $G\oplus H$. In addition, if $a(G)=1$, then $S_G\oplus (h)$ is a double dominating sequence of $G\oplus H$. This fact proves that $\dgri{G\oplus H}\geq\dgri{G}+a(G)$.\\
    
    Let us show that $\dgri{G\oplus H}\leq\max\{\dgri{G}+a(G),\dgri{H}+a(H)\}$. Consider $S$ a MDNS of $G\oplus H$, $S_G=S-V(H)=(g_1,\dots,g_k)$, and $S_H=S-V(G)=(h_1,\dots,h_{m})$. Without loss of generality, assume that $|S_G|\geq |S_H|$.
    Consider $|\widehat{S_G}|=1$. If $a(G)=1$ or $|\widehat{S_G}|\leq\dgri{G}-1$, the result is trivial. Assume that $a(G)=0$ and $S_G$ is a MDNS of $G$. Observe that $|\widehat{S_G}|\geq2$. In addition, $|N_{G\vee H}[v]\cap\widehat{S}|\geq2$ for all $v\in V(G\vee H)$, then $o_S(g_k)>o_S(h_1)$. If $k\geq3$, we have $|N_{G\vee H}[v]\cap\{s_1,\dots,s_{k-1},h_1\}|\geq2$ for all $v\in V(G\vee H)$ and then $N_S[g_k]=\emptyset$, which contradicts that $S$ is a MDNS of $G\vee H$. If $k=2$, then $\dgri{G}=2$ and $|\widehat{S}|=3$. Since $a(G)=0$, it is not hard to see that $G$ is a complete graph. Hence, $o_S(g_2)=3$. Consequently, $N_S[g_2]=N_S^2[g_2]\subset V(H)$ and then there exists $h\in V(H)-N_H[h_1]$. Hence, $\dgri{H}\geq2$. If $a(H)=1$, $|\widehat{S}|=3\leq\dgri{H}+a(H)$. If $a(H)=0$, let $h'\in N_H[h]$. Thus, $(h',h_1,h)$ is a DNS of $H$ and $|\widehat{S}|=3\dgri{H}+a(H)$.
    %and thus $G$ then $o_S(g_k)>o_S(h)$. $o_S(h_1)>o_S(g_2)$, then $N_S[h_1]\subseteq V(G)$. If  Hence, $o$ $a(G)=1$ or if $a(G)=0$, then $S_G$ is not a MDNS of $G$, that is, $|\widehat{S_G}|\leq\dgri{G}-1$. Thus, $|\widehat{S}|\leq\dgri{G}+a(G)$. If $o_S(h_1)<o_S(g_2)$, since $V(G)\subseteq N_{G\vee H}[h_1]$, then $V(G)\subseteq N_S[h_1]$. As before, $a(G)=1$ or if $a(G)=0$, then $S_G$ is not a MDNS of $G$. Thus, $|\widehat{S}|\leq\dgri{G}+a(G)$.

    Consider now $|\widehat{S_H}|=2$. Then, $|\widehat{S_G}|\geq2$. If $o_S(g_2)<o_S(h_2)$, $N_S[h_2]=N_S^2[h_2]\subseteq V(G)$. Let $v\in N_S^2[h_2]$. Thus, $v\in N_S^1[h_1]$. Hence, if $v$ is not an isolated vertex of $G$ and $v'\in N_G(v)$, then the sequence $S'$ obtained replacing $h_1$ and $h_2$ in $S$ by $v$ and $v'$ is a DNS of $G$. Then, $|\widehat{S_G}|\leq\dgri{G}-2$. Otherwise, if $v$ is an isolated of $G$, $v\notin\widehat{S_G}$, $|\widehat{S_G}|\leq\dgri{G}-1$, and $a(G)=1$. Thus, in both cases, $|\widehat{S}|\leq\dgri{G}+a(G)$. If $o_S(g_2)>o_S(h_2)$, then $g_2$ is the last vertex of $S$, since for $i=1,2$, $V(G)\subseteq N_{G\vee H}[h_i]$ and $V(H)\subseteq N_{G\vee H}[g_i]$. In consequence, $|\widehat{S}|=4$. In addition, analogously as before, $N_S[g_2]=N_S^2[g_2]\subseteq V(H)$. Let $h\in N_S^2[g_2]$. Thus, $h\in N_S^1[g_1]$. If $h$ is not an isolated vertex of $H$ then $|\widehat{S_H}|=2\leq\dgri{H}-2$ and, if $h$ is an isolated of $H$, then $h\notin\widehat{S_H}$, $|S_H|=2\leq\dgri{H}-1$, and $a(H)=1$. Thus, $|\widehat{S}|\leq\dgri{H}+a(H)$.

    Finally, assume that $|S_G|\geq3$. Then, $N_S[g_3]\subseteq V(G)$. If $|S_H|\geq3$, $N_S[h_3]\subseteq V[H]$. If $o_S(h_3)<o_S(g_3)$, since $V(G)\subset N_{G\vee H}[h_1]$ and $V(G)\subset N_{G\vee H}[h_2]$, $N_S[g_3]=\emptyset$, which is a contradiction. Analogously if $o_S[h_3]>o_S[g_3]$. Therefore, $|S_H|\leq2$ and the result is derived from the cases analyzed above.

    In addition, it is not hard to see that, in each case, the sequences given in items $(a)-(d)$ are MDNSs of $G\vee H$.
\end{proof}

Given a graph $G$, it is immediate to see that if $G$ (resp. $\overline G$) is not connected, $G$ can be obtained by disjoint union (resp. join) of two non-empty induced subgraphs of $G$. When $G$ and $\overline{G}$ are connected, $G$ is called \emph{modular}. Then, given a graph family $\mathcal F$, $M(\mathcal F)$ denotes the family of modular graphs in $\mathcal F$.

Here we consider the modular decomposition following \cite{BKNT-2018}. Notice that Lemma \ref{lemma:unionyjoin} implies that a MDNS of $G_1+ G_2$ and $G_1\vee G_2$ can be obtained in linear time from MDNSs of $G_1$ and $G_2$. This implies the following.

\begin{rem}\label{rem:modular}{\cite{BKNT-2018}}
    Given a graph class $\mathcal F$, if a MDNS can be computed in linear (polynomial) time for every graph in $M(\mathcal F)$, then the problem of finding a MDNS is linear (polynomial) time for graphs in $\mathcal F$.
\end{rem}

In this regard, it is known that if $\mathcal F$ is the class of cographs, then every graph in $M(\mathcal F)$ is trivial. Since for a trivial graph $G=(\{v\},\emptyset)$ the only MDNS is $S=(v)$, we see that the problem of finding a MDNS for cographs is linear. This implies the following result.

\begin{theo}\label{theo:cographs}
    The GDDN and a GDDS can be obtained in linear time for cographs. 
\end{theo}

Cographs serve as the base of a hierarchy of graph classes defined by their $P_4$ structure, being precisely the $P_4$-free graphs. Superclasses of cographs include extended $P_4$-sparse~\cite{GV-1997}, $P_4$-lite~\cite{JO3}, and $P_4$-extendible graphs~\cite{JO4}. These are all generalized by the class of $P_4$-tidy graphs~\cite{Giako}. The next section will explore Grundy double dominating sequences in this family.

%%%%%%%%%%%%%%%%%%%%%%%%%%%%%%%%%%%%%%%%%%%%%%%%%%%%%%%%%%%%%%%%%%%%%%%%
\subsection{$P_4$-tidy graphs}\label{sec:tidy}

Let us introduce a graph class that play an important role in this section. A graph $G$ is considered a \emph{spider graph} if its vertices can be partitioned into three sets, $S$, $C$, and $H$ (where $H$ can be empty), satisfying the following properties: $S$ is a stable set, $C$ is a clique, $|S|=|C|=r\geq2$, every vertex in $C$ is adjacent to all vertices in $H$, and no vertex in $H$ is adjacent to any vertex in $S$. Furthermore, if $S=\{s_1, \ldots, s_r\}$ and $C= \{c_1, \ldots, c_r\}$, one of these conditions must hold:
\begin{enumerate}
\item \emph{Thin spider}: $s_i$ is adjacent to $c_j$ if and only if $i = j$.
\item \emph{Thick spider}: $s_i$ is adjacent to $c_j$ if and only if $i \neq j$.
\end{enumerate}

The cardinality of $C$ (and $S$) is referred to as the \emph{weight} of $G$, and the set $H$ is known as the \emph{head} of the spider. A spider graph $G$ with the partition $S,C,H$ is denoted as $G=(S,C,H)$.

It is worth noting that for $r=2$, both thin and thick spider graphs are isomorphic to the path $P_4$. Therefore, from this point forward, we will exclusively consider thick spider graphs with $r\geq3$.

Let $G=(S,C,H)$ a thin (resp. thick) spider graph. The graph obtained by adding a true twin or a false twin vertex of a vertex $v\in S\cup C$ is called \emph{thin (resp. thick) quasi-spider graph}. Without loss of generality, we assume that $v$ is $s_r$ or $c_r$ and the added vertex is called $s'_r$ or $c'_r$, respectively.
In addition, we denote by $(S_f,C,H)$ and $(S_t,C,H)$ the quasi-spider graph obtained from $G=(S,C,H)$ by adding to $s_r$ the false or true twin $s'_r$, respectively. Analogously, $(S,C_f,H)$ and $(S,C_t,H)$ are the quasi-spider graph obtained from $G=(S,C,H)$ by adding to $c_r$ the false or true twin $c'_r$, respectively. In addition, we consider the weight of a quasi-spider graph as the weight of the original spider graph, $S_f=S_t=S\cup\{s'_r\}$, and $C_f=C_t=C\cup\{c'_r\}$. 

For spider and quasi-spider graphs, their partition is uniquely determined, and both their recognition and the determination of their partition can be achieved in linear time~\cite{Giako}.

%From Proposition \ref{cor:twins} we infer that $\grt (S\hookleftarrow \overline K_2,C,H)=\grt (S,C\hookleftarrow \overline K_2,H)=\grt (S,C,H)$. Thus it only remains to compute the Grundy total domination number for spider graphs and quasi-spider graphs of the type $(S\hookleftarrow K_2,C,H)$ and $(S,C\hookleftarrow K_2,H)$.

In order to determine the computational complexity of computing the GDDN of a $P_4$-tidy graph, we applied an approach similar to that used in the previous section according to Remark~\ref{rem:modular}. Hence, we are interested in obtaining GDNSs of modular $P_4$-tidy graphs.  In \cite{Giako} the family of modular $P_4$-tidy is determined: modular $P_4$-tidy graphs are trivial graphs, $P_5$, $\overline{P_5}$ (house), $C_5$, and spider graphs and quasi-spider graphs such that their head is a $P_4$-tidy graph.

In the following results, we determine the maximum double neighborhood number of spider and quasi-spider graphs. Recall that it coincides with the GDDN if the graph does not have isolated vertices. To simplify, if $H=\emptyset$ we consider $S_H=()$ a MDNS of $G[H]$ and $\dgri{G[H]}=0$.

\begin{lemma}\label{lemma:gdriH=0}
    Let $G=(S,C,\emptyset)$ be a spider graph with weight $r$. If $G$ is a thin spider graph, then $\dgri{G}=2r$ and $(c_1,\dots,c_r,s_1,\dots,s_r)$ is a MDNS of $G$. If $G$ is a thick spider graph, then $\dgri{G}=2+r$ and $(s_1,\dots,s_r,c_1,c_2)$ is a MDNS of $G$.
\end{lemma}
\begin{proof}
    If $G$ is a thin spider graph, the result follows from the fact that $(c_1,\dots,c_r,s_1,\dots,s_r)$ is a MDNS of $G$.

    Assume now that $G$ is a thick spider graph. It is easy to see that $(s_1,\dots,s_r,c_1,c_2)$ is a DNS of $G$. Let $T$ be a MDNS of $G$ and $T_C=T-S=(c_{i_1},c_{i_2},\dots,c_{i_k})$. If $k\leq2$ it is immediate that $|\widehat{T}|\leq2+r$. If $k\geq3$, then $N_T[c_{i_3}]=N_T^2[c_{i_3}]\subseteq\{s_{i_1},s_{i_2}\}$. In consequence $k=3$. In addition, observe that for all $v\in V(G)$, $|N_G[v]\cap \{c_{i_1},c_{i_2},c_{i_3}\}|\geq2$. Hence, $c_{i_3}$ is the last vertex in the sequence $T$, $\{s_{i_1},s_{i_2}\}\nsubseteq\widehat{T} $, and then $|\widehat{T}\cap S|\leq r-1$. Thus, $|\widehat{T}|\leq2+r$. Therefore, $(s_1,\dots,s_r,c_1,c_2)$ is a MDNS of $G$ and $\dgri{G}=2+r$.
\end{proof}

\begin{lemma}\label{lemma:T_H}
    Let $G=(S,C,H)$ be a spider graph and $T$ a MDNS of $G$. Then $|\widehat{T}\cap H|\leq\dgri{G[H]}$ and $|\widehat{T}\cap (C\cup S)|\leq\dgri{G[C\cup S]}$.
\end{lemma}
\begin{proof}
    Let $T$ be a MDNS of $G$.    
    Consider $T_H=T-(C\cup S)$. Hence, $\widehat{T_H}=\widehat{T}\cap H$. Let us show that $|\widehat{T_H}|\leq\dgri{G[H]}$. If $|\widehat{T_H}|\leq2$, the result follows immediately. Assume that $k=|\widehat{T_H}|\geq3$ and let $T_H=(v_1,v_2,\dots,v_k)$. Since $C\subset N_G[v]$ for all $v\in H$, then $N_T[v_i]\subseteq N_{G[H]}[v_i]$ for all $i=3,\dots,k$. Therefore, $T_H$ is a DNS of $G[H]$ and $|\widehat{T}\cap H|=|\widehat{T_H}|\leq\dgri{G[H]}$.

    Observe that $G[C\cup S]=(S,C,\emptyset)$.
    If $G$ is a thin spider graph, the result follows immediately from Lemma~\ref{lemma:gdriH=0}.

    Finally, consider $G$ a thick spider graph and $T_C=T-(H\cup S)=(c_{i_1},c_{i_2},\dots,c_{i_k})$. If $k\leq2$, the result follows from Lemma~\ref{lemma:gdriH=0}. If $k\geq3$, then $N_T[c_{i_3}]=N_T^2[c_{i_3}]\subseteq\{s_{i_1},s_{i_2}\}$. In consequence, $k=3$ and by a reasoning similar to that in Lemma~\ref{lemma:gdriH=0}, $|T|\leq2+r=\dgri{G}$.
\end{proof}

From Lemma~\ref{lemma:T_H}, we have $\dgri{G}\leq\dgri{G[C\cup S]}+\dgri{G[H]}$. Let $T_H$ be a MDNS of $G[H]$. If $G$ is a thin spider graph, then $T_H\oplus (c_1,\dots,c_r,s_1,\dots,s_r)$ is a MDNS of $G$. Similarly, if $G$ is a thick spider graph, $T_H\oplus (s_1,\dots,s_r,c_1,c_2)$ is a MDNS of $G$ and the following is immediate.

\begin{prop}\label{prop:spider}
    Let $G=(S,C,H)$ be a spider graph with weight $r$ and $T_H$ be a MDNS of $G[H]$. Then, \[\dgri{G}=\dgri{G[C\cup S]}+\dgri{G[H]}.\]    
    In addition, if $G$ is thin, $T_H\oplus (c_1,\dots,c_r,s_1,\dots,s_r)$ is a MDNS of $G$ and if $G$ is thick, $T_H\oplus (s_1,\dots,s_r,c_1,c_2)$ is a MDNS of $G$.
\end{prop}

In what follows, we analize MDNSs of quasi-spider graphs. To this end, we provide two useful results.

\begin{lemma}\label{rem:CbeforeH}
    Let $G=(S,C,H)$ be a spider graph such that $H\neq\emptyset$ and $a(G[H])=0$. Consider $T$ a DNS of $G$ such that there exist $c\in C\cap\widehat{T}$ and $h\in H\cap\widehat{T}$ with $o_T(c)<o_T(h)$. Then, $|H\cap\widehat{T}|\leq\dgri{G[H]}-1$. 
\end{lemma}
\begin{proof}
    Let $T_H=T-(C\cup S)=(h_1,\dots,h_k)$. Note that $\widehat{T_H}=H\cap\widehat{T}$. Since $H\neq\emptyset$ and $a(G[H])=0$, then $|\widehat{T_H}|\geq2$. Hence, $o_T(h_k)\geq o_T(h)>o_T(c)$. Thus, $N_T[h_k]=N_T^2[h_k]\subseteq H$. Let $h'\in N_T^2[h_k]$. It follows that $h'\in N_T^1[c]$. Consequently, $h'=h_k$ or $h'\notin\widehat{T_H}$. If $h'=h_k$, since $a(H)=0$, there exists $h''\in N_{G[H]}(h')$. Then, $T_H'=T_H\oplus (h'')$ is a DNS of $G[H]$, since $h_k\in N_{T_H'}[h'']$. If $h'\notin\widehat{T_H}$, $T_H'=T_H\oplus (h')$ is a DNS of $G[H]$, since $h_k\in N_{T_H'}[h']$. Therefore, in both cases, $|\widehat{T_H}|=|\widehat{T'_H}|-1\leq\dgri{G[H]}-1$ and the result holds.
\end{proof}

A similar reasoning as in Lemma~\ref{lemma:T_H} for quasi-spider graph implies the following.

\begin{rem}\label{rem:T_Hquasi}
    Let $G=(S,C,H)$ be a spider graph, $G'$ a quasi-spider graph obtained from $G$, and $T$ a MDNS of $G'$. Then, $|\widehat{T}\cap H|\leq\dgri{G[H]}$.
\end{rem}

\begin{prop}\label{prop:thinquasispider}
    Let $G=(S,C,H)$ be a thin spider graph with weight $r$ and $T_H$ be a MDNS of $G[H]$. Then:
        \begin{enumerate}
            \item $\dgri{(S_f,C,H)}=\dgri{G}+1=2r+1+\dgri{G[H]}$ and $T=T_H\oplus (c_1,\dots,c_r,s_1,\dots,s_r,s'_r)$ is a MDNS of $(S_f,C,H)$.
            \item If $H=\emptyset$ or $a(G[H])=1$, then $\dgri{(S_t,C,H)}=\dgri{G}+1+\dgri{G[H]}=2r+1+\dgri{G[H]}$ and $T=T_H\oplus (s_r,s'_r,c_r,c_1,\dots,c_{r-1},s_1,\dots,s_{r-1})$ is a MDNS of $(S_t,C,H)$.
            \item If $H\neq\emptyset$ and $a(G[H])=0$, then $\dgri{(S_t,C,H)}=\dgri{G}+\dgri{G[H]}=2r+\dgri{G[H]}$ and $T=T_H\oplus (c_1,\dots,c_r,s_1,\dots,s_r)$ is a MDNS of $(S_t,C,H)$.
            \item If $H=\emptyset$, then $\dgri{(S,C_f,H)}=\dgri{(S,C_t,H)}=\dgri{G}+1=2r+1$ and $T=(s_r,c_r,c'_r,c_1,\dots,c_{r-1},s_1,\dots,s_{r-1})$ is a MDNS of $(S,C_f,H)$ and $(S,C_t,H)$.
            \item If $H\neq\emptyset$, then $\dgri{(S,C_f,H)}=\dgri{(S,C_t,H)}=\dgri{G}+\dgri{G[H]}=2r+\dgri{G[H]}$ and $T=T_H\oplus (c_1,\dots,c_r,s_1,\dots,s_r)$ is a MDNS of $(S,C_f,H)$ and $(S,C_t,H)$.
    \end{enumerate}
\end{prop}
\begin{proof}
    \begin{enumerate}
        \item From propositions~\ref{prop:false} and \ref{prop:spider}, $\dgri{(S_f,C,H)}\leq\dgri{G}+1=2r+1+\dgri{G[H]}$. Furthermore, it is not hard to see that $T=T_H\oplus (c_1,\dots,c_r,s_1,\dots,s_r,s'_r)$ is a DNS of $(S_f,C,H)$. Thus, $2r+1+\dgri{G[H]}=|\widehat{T}|\leq\dgri{(S_f,C,H)}\leq2r+1+\dgri{G[H]}$ and the result follows.
        
        \item Analogously as above, from propositions~\ref{prop:true} and \ref{prop:spider}, $\dgri{(S_t,C,H)}\leq\dgri{G}+1=2r+1+\dgri{G[H]}$. Furthermore, it is not hard to see that $T=T_H\oplus (s_r,s'_r,c_r,c_1,\dots,c_{r-1},s_1,\dots,s_{r-1})$ is a DNS of $(S_t,C,H)$. Thus, $2r+1+\dgri{G[H]}=|\widehat{T}|\leq\dgri{(S_t,C,H)}\leq2r+1+\dgri{G[H]}$ and the result follows.
        
        \item It is easy to see that $T=T_H\oplus (c_1,\dots,c_r,s_1,\dots,s_r)$ is a DNS of $(S_t,C,H)$. If $T'$ is a MDNS of $(S_t,C,H)$ with $|\widehat{T'}|>|\widehat{T}|=2r+\dgri{G[H]}$, from Proposition~\ref{prop:true}, $|\widehat{T'}|=2r+1+\dgri{G[H]}$ and there exists a MDNS $T''$ of $G$ such that $s_r\in\widehat{T''}$, $N_{T''}^1[s_r]\neq\emptyset$ (*), and for every $x\in\widehat{T''}-\{s_r\}$, $N_{T''}[x]-N^1_{T''}[s_r]\neq\emptyset$ (**). Furthermore, from Lemma~\ref{lemma:T_H}, $T''_H=T''-(C\cup S)$ is a MDNS of $G[H]$. Therefore, since $c_r\in\widehat{T''}$ and $N_G[s_r]=\{s_r,c_r\}\subset N_G[c_r]$, from (*) we have $o_{T''}(s_r)<o_{T''}(c_r)$. Since $H\neq\emptyset$ and $a(H)=0$, then $|\widehat{T''_H}|=\dgri{G[H]}\geq2$ and, by Lemma~\ref{rem:CbeforeH}, $o_{T''}(h)<o_{T''}(c_r)$ for every $h\in\widehat{T''_H}$. Thus, $N_{T''}[c_r]=\{s_r\}$, which contradicts (**). Hence, such a sequence $T'$ does not exist and $T$ is a MDNS of $(S_t,C,H)$.
        \item The result is trivial since $T=(s_r,c_r,c'_r,c_1,\dots,c_{r-1},s_1,\dots,s_{r-1})$ verifies that $\widehat{T}=S\cup C\cup\{c'_r\}$ and it is a DNS of $(S,C_f,H)$ and $(S,C_t,H)$. Thus, $\dgri{(S,C_f,H)}=\dgri{(S,C_t,H)}=\dgri{G}+1=2r+1$ and $T$ is a MDNS of $(S,C_f,H)$ and $(S,C_t,H)$.
        
        \item It is easy to see that $T=T_H\oplus (c_1,\dots,c_r,s_1,\dots,s_r)$ is a DNS of $(S,C_f,H)$ and $(S,C_t,H)$. Let $G'\in \{(S,C_f,H),(S,C_t,H)\}$. If there exists a MDNS $T'$ of $G'$ with $|\widehat{T'}|>|\widehat{T}|$, from propositions~\ref{prop:true} and \ref{prop:false}, $|\widehat{T'}|=2r+1+\dgri{G[H]}$. By Remark~\ref{rem:T_Hquasi}, $T'_H=T'-(C\cup S)$ is a MDNS of $G[H]$ and $C\cup S\subset\widehat{T'}$. Note that $\dgri{G[H]}\geq1$. Consider the vertices $s_r,c_r,c'_r,h$ with $h\in\widehat{T'_H}$. We have that the four vertices belong to $\widehat{T'}$. In addition, $N_{G'}[s_r]=\{s_r,c_r,c'_r\}$, $(C_t-\{c'_r\})\cup H\cup\{s_r\}\subseteq N_{G'}[c_r]$, $(C_t-\{c_r\})\cup H\cup\{s_r\}\subseteq N_{G'}[c'_r]$, and $C_t\cup N_{G[H]}[h]=N_{G'}[h]$. Therefore, among $s_r,c_r,c'_r$, and $h$, either $c_r$ or $c'_r$ has the largest relative order in $T'$. Without loss of generality, we assume that $o_{T'}(h),o_{T'}(s_r),o_{T'}(c_r)<o_{T'}(c'_r)$ for every $h\in\widehat{T'}\cap H$. Hence, $N_{T'}[c_r']=N_{T'}^2[c'_r]\subseteq H$. Let $w\in N_{T'}^2[c'_r]$. Since $w\in N_{G'}[c_r]$, then $w\in N^1_{T'}[c_r]$. These facts imply that $N_{G[H]}[w]\cap\widehat{T'}=\emptyset$, which contradicts the fact that $T'_H$ is a MDNS of $G[H]$.
        Thus, $T=T_H\oplus (c_1,\dots,c_r,s_1,\dots,s_r)$ is a MDNS of $(S,C_f,H)$ and $(S,C_t,H)$.
    \end{enumerate}
\end{proof}

\begin{prop}\label{prop:thickquasispider}
    Let $G=(S,C,H)$ be a thick spider graph with weight $r$ and $T_H$ is a MDNS of $G[H]$. Then:
    \begin{enumerate}
        \item $\dgri{(S_f,C,H)}=\dgri{(S_t,C,H)}=\dgri{G}+1=r+3+\dgri{G[H]}$ and $T=T_H\oplus (s_1,\dots,s_r,s'_r,c_1,c_2)$ is a MDNS of $(S_f,C,H)$ and $(S_t,C,H)$.
        \item $\dgri{(S,C_f,H)}=\dgri{(S,C_t,H)}=\dgri{G}=r+2+\dgri{G[H]}$ and $T=T_H\oplus (s_1,\dots,s_r,c_1,c_2)$ is a MDNS of $(S,C_f,H)$ and $(S,C_t,H)$.
    \end{enumerate}
\end{prop}
\begin{proof}
    \begin{enumerate}
        \item It is easy to see that $T=T_H\oplus (s_1,\dots,s_r,s'_r,c_1,c_2)$ is a DNS of $(S_f,C,H)$ and $(S_t,C,H)$. The result is derived from propositions \ref{prop:true}, \ref{prop:false}, and \ref{prop:spider}.
        
        \item We first analyze $G'=(S,C_t,H)$. It is easy to see that $T=T_H\oplus (s_1,\dots,s_r,c_1,c_2)$ is a DNS of $(S,C_t,H)$. If $T'$ is a MDNS of $(S,C_t,H)$ with $|\widehat{T'}|>|\widehat{T}|=r+2+\dgri{G[H]}$, from Proposition~\ref{prop:true}, $|\widehat{T'}|=r+3+\dgri{G[H]}$ and there exists a MDNS $T''$ of $G$ such that  $|\widehat{T''}|=r+2+\dgri{G[H]}$, $c_r\in\widehat{T''}$, $N_{T''}^1[c_r]\neq\emptyset$ (*), and for every $x\in\widehat{T''}-\{c_r\}$, $N_{T''}[x]-N^1_{T''}[c_r]\neq\emptyset$ (**). 
        Let $T''_C=T''-(S\cup H)$. Using a reasoning similar to that in the proof of Lemma \ref{lemma:gdriH=0}, it follows that $|\widehat{T''_C}|\in\{2,3\}$.
        
        If $|\widehat{T''_C}|=2$, then $S\subset\widehat{T''}$ and $|\widehat{T''}\cap H|=\dgri{G[H]}$. From (**), $o_{T''}(s_i)<o_{T''}(c_r)$ for all $i\in\{1,\dots,r-1\}$ ($r\geq3$). Then, $N^1_{T''}[c_r]\subseteq H$. Hence, $H\neq\emptyset$. Since $N_G[h]\subset N_G[c_r]$, from (**), $o_{T''}(h)<o_{T''}(c_r)$ for all $h\in\widehat{T''}\cap H$.  In addition, $|\widehat{T''}\cap H|=\dgri{G[H]}$. Then, $N^1_{T''}(c_r)=\emptyset$, which contradicts (*).
            
        Assume now that $|\widehat{T''_C}|=3$ and $T''_C=(w_1,w_2,w_3)$. Using a reasoning similar to that in the proof of Lemma \ref{lemma:T_H}, it follows that $|\widehat{T''}\cap S|=r-1$ and $|\widehat{T''}\cap H|=\dgri{G[H]}$. If $\{s_1,\dots,s_{r-1}\}=\widehat{T''}\cap S$, from (**), $o_{T''}(s_i)<o_{T''}(c_r)$ for all $i\in\{1,\dots,r-1\}$. Following a similar argument as above, we conclude the thesis. Then, without loss of generality, we assume that $\{s_1,\dots,s_{r-2},s_r\}=\widehat{T''}\cap S$. Hence, $o_{T''}(s_i)<o_{T''}(c_r)$ for all $i\in\{1,\dots,r-2\}$. If $H=\emptyset$, $N^1_{T''}[c_r]=\{s_{r-1}\}$. Otherwise, again from (**), $o_{T''}(h)<o_{T''}(c_r)$ for all $h\in\widehat{T''}\cap H$.  In addition, $|\widehat{T''}\cap H|=\dgri{G[H]}$. Thus, in any case, $N_{T''}^1[c_r]=\{s_{r-1}\}$. Consequently, $c_r\in\{w_1,w_2\}$. We also know that $o_{T''}(s_r)<o_{T''}(w_3)$. From the fact that $N_G[s_r]\cap\{w_1,w_2\}\neq\emptyset$, we have $N_{T''}[w_3]=\{s_{r-1}\}= N^1_{T''}[c_r]$, which contradicts (**).
        
        Therefore, in both cases, such a sequence $T'$ does not exist and $T$ is a MDNS of $(S,C_t,H)$.
%        Using a reasoning similar to that in the proof of Lemma \ref{lemma:gdriH=0}, it follows that $|\widehat{T''}\cap C|=3$ and $|\widehat{T''}\cap S|=r-1$. Let $T''-(C\cup S)=(w_1,w_2,w_3)$. By (*), it is clear that $c_r\in\{w_1,w_2\}$. Since for $i\in\{1,\dots,r-1\}$, $N_G[s_i]\subset N_G[c_r]$, from (**) follows that each vertex in $\widehat{T''}\cap\{s_1,\dots,s_{r-1}\}$ appears in $T''$ before $c_r$. Note that $N_G[c_r]=\cup_{i=1}^{r-1}N_G[s_i]$ and then, by condition (*), $\widehat{T''}\cap S\neq\{s_1,\dots,s_{r-1}\}$, that is, $s_r\in\widehat{T''}\cap S$. W.l.o.g. assume that $\widehat{T''}\cap S=\{s_1,\dots,s_{r-2},s_r\}$. Hence, $o_{T''}(s_i)<o_{T''}(c_r)$ for every $i\in\{1,\dots,r-2\}$. Thus, $N_{T''}^1(c_r)=\{s_{r-1}\}$ and $o_{T''}(s_r)<o_{T''}(w_3)$. From the fact that $s_r\in\widehat{T''}$ and $N_G[s_r]\cap\{c_{i_1},c_{i_2}\}$, then $N_{T''}[c_{i_3}]=\{s_{r-1}\}\subset N_{T''}[c_r]$, which contradicts (**). 

       Consider now $G'=(S,C_f,H)$. It is easy to see that $T=T_H\oplus (s_1,\dots,s_r,c_1,c_2)$ is a DNS of $G'$. If $T'$ is a MDNS of $G'$ with $|\widehat{T'}|>|\widehat{T}|=r+2+\dgri{G[H]}$, from Proposition~\ref{prop:false} and Remark \ref{rem:T_Hquasi}, $|\widehat{T'}|=r+3+\dgri{G[H]}$, $|\{c_r,c'_r\}\cap\widehat{T'}|\geq1$, and $|\widehat{T'}\cap H|\leq\dgri{G[H]}$. Without loss of generality, we assume that $c_r\in\widehat{T'}$. Let $T'_C=T'-(S\cup H)=(w_1,\dots,w_k)$. Observe that $k\in\{3,4\}$, since $|N_{G'}[v]\cap\{w_1,w_2,w_3,w_4\}|\geq2$, for all $v\in V(G')$. We consider the cases $c'_r\notin\widehat{T'}$, $c'_r\in\widehat{T'}$ and $k=3$, and $c'_r\in\widehat{T'}$ and $k=4$.

        If $c'_r\notin\widehat{T'}$, then $|N_{G'}[v]\cap\{w_1,w_2,w_3\}|\geq2$ for all $v\in V(G')$. We have $k=3$, $S\subset\widehat{T'}$ and $w_3$ is the last vertex of $T'$. Furthermore, for all $v\in V(G')$, $|N_{G'}[v]\cap(S\cup\{w_1,w_2\})|\geq2$. Thus, $N_{T'}[w_3]=\emptyset$, which contradicts the fact that $T'$ is a DNS of $G'$. Therefore, $T$ is a MDNS of $(S,C_f,H)$.

        Consider now $c'_r\in\widehat{T'}$. Without loss of generality, we assume that $o_{T'}(c_r)<o_{T'}(c'_r)$.% In addition, $|\widehat{T'}\cap C_f|\in\{3,4\}$ and $|\widehat{T'}\cap S|\in\{r,r-1\}$. 
        
        If $k=3$, then $|\widehat{T'}\cap H|=\dgri{G[H]}$ and $S\subset\widehat{T'}$. Among the vertices $s_1,\dots,s_{r-1}$, let $s$ be the largest relative order vertex in $T'$. Since $r\geq3$, there exists $s'\in\{s_1,\dots,s_{r-1}\}-\{s\}$ with $o_{T'}(s')<o_{T'}(s)$. Observe that every vertex in $N_{G'}[s']$ is dominated at least twice by $\{s,c_r,c'_r\}$. Consequently, for every $i=1,\dots,r-1$, $o_{T'}(s_i)<o_{T'}(c'_r)$ ($r\geq3$). Therefore, $N_{T'}[c'_r]=N^2_{T'}[c'_r]\subseteq H$. Hence, there exists $h\in N^2_{T'}[c_r']\subseteq H$. Since $|N_{G'}[v]\cap (\{s_1,\dots,s_{r-1}\}\cup\{c_r,c'_r\}|\geq2$ for all $v\in V(G')-\{s_r\}$, then $o_{T'}(h')<o_{T'}(c'_r)$ for all $h'\in\widehat{T'}\cap H$. But $h\in N_{G'}[c_r]$ and $o_{T'}(c_r)<o_{T'}(c'_r)$, thus $N_{G[H]}[h]\cap\widehat{T'}=\emptyset$. Hence, $|\widehat{T'}\cap H|<\dgri{G[H]}$, which is a contradiction. Thus, $T$ is a MDNS of $(S,C_f,H)$.
        
        Finally, let $k=4$. Recall that $\{c_r,c'_r\}\subset\widehat{T'}$ and $o_{T'}(c_r)<o_{T'}(c'_r)$. Observe that $o_{T'}(c'_r)<o_{T'}(w_4)$ and $N_{T'}[w_4]=N_{T'}^2[w_4]=\{s_r\}$. Hence, if $s_r\in\widehat{T'}$, since $|N_{G'}[s_r]\cap\{w_1,w_2,w_3\}|=1$, then $o_{T'}(w_4)<o_{T'}(s_r)$. But $|N_{G'}[s_r]\cap\{w_1,w_2,w_3,w_4\}|=2$, then $s_r\notin\widehat{T'}$. Therefore, $\{s_1,\dots,s_{r-1}\}=\widehat{T'}\cap S$ and so, $|\widehat{T'}\cap H|=\dgri{G[H]}$. Thus, as we have previously observed, for every $i=1,\dots,r-1$, $o_{T'}(s_i)<o_{T'}(c'_r)$ ($r\geq3$), $N_{T'}[c'_r]=N^2_{T'}[c'_r]\subseteq H$, and we conclude that $|\widehat{T'}\cap H|<\dgri{G[H]}$, which is a contradiction. Thus, $T$ is a MDNS of $(S,C_f,H)$.
        %for all $i=1,\dots,r-1$, $o_{T'}(s_i)<o_{T'}(c'_r)$ and $N_{T'}[c'_r]=N^2_{T'}[c'_r]\subseteq H$. Let $h\in N_{T'}[c'_r]$ If $H=\emptyset$ the result follows. Else, since $H\subseteq N[c_r]$, $|\widehat{T'}\cap H|<\dgri{G[H]}$ and then such a set $T'$ does not exist. Therefore, $T$ is a MDNS of $(S,C_f,H)$.
    \end{enumerate}
\end{proof}

Given that non-trivial modular $P_4$-tidy graphs consist of spider and quasi-spider graphs with $P_4$-tidy heads, along with $C_5$, $P_5$ and $\bar P_5$, we can directly derive the following theorem by combining Remark~\ref{rem:modular}, Propositions~\ref{prop:spider}, \ref{prop:thinquasispider} and \ref{prop:thickquasispider}.

\begin{theo}\label{theo:P4tidy}
The GDDN and a GDDS can be obtained in linear time for $P_4$-tidy graphs.
\end{theo}

\section*{Declarations}
\textbf{Conflict of interest:} The author has no conflicts of interest to declare that are relevant to the content of this article.

\bibliographystyle{apalike}
\bibliography{bibliography}

\begin{thebibliography}{}

\bibitem[Bre\v{s}ar et~al., 2014]{BGMRR-2014}
Bre\v{s}ar, B., Gologranc, T., Milani\v{c}, M., Rall, D.~F., and Rizzi, R. (2014).
\newblock Dominating sequences in graphs.
\newblock {\em Discrete Mathematics}, 336:22--36.

\bibitem[Bre\v{s}ar et~al., 2010]{BKR-2010}
Bre\v{s}ar, B., Klav{\v{z}}ar, S., and Rall, D.~F. (2010).
\newblock Domination game and an imagination strategy.
\newblock {\em SIAM Journal of Discrete Mathematics}, 24:979--991.

\bibitem[Bre\v{s}ar et~al., 2023]{BSP-2023}
Bre\v{s}ar, B., Pandey, A., and Sharma, G. (2023).
\newblock Computation of grundy dominating sequences in (co-)bipartite graphs.
\newblock {\em Comput. Appl. Math.}, 42(8): 359.

\bibitem[Brešar et~al., 2021]{BHKR-2021}
Brešar, B., Henning, M.~A., Klavžar, S., and Rall, D.~F. (2021).
\newblock {\em Domination Games Played on Graphs}.
\newblock Springer Cham, Switzerland.

\bibitem[Brešar et~al., 2018]{BKNT-2018}
Brešar, B., Kos, T., Nasini, G., and Torres, P. (2018).
\newblock Total dominating sequences in trees, split graphs, and under modular decomposition.
\newblock {\em Discrete Optimization}, 28:16--30.

\bibitem[Brešar et~al., 2022]{BPS-2021}
Brešar, B., Pandey, A., and Sharma, G. (2022).
\newblock Computational aspects of some vertex sequences of grundy domination-type.
\newblock {\em Indian J. Discrete Math.}, 8:21--38.

\bibitem[Cournier and Habib, 1994]{CH-1994}
Cournier, A. and Habib, M. (1994).
\newblock A new linear algorithm for modular decomposition.
\newblock In Tison, S., editor, {\em Trees in Algebra and Programming --- CAAP'94}, pages 68--84, Berlin, Heidelberg. Springer Berlin Heidelberg.

\bibitem[Giakoumakis et~al., 1997]{Giako}
Giakoumakis, V., Roussel, F., and Thuillier, H. (1997).
\newblock On $p_4$-tidy graphs.
\newblock {\em Discrete Mathematics and Theoretical Computer Science}, 1:17--41.

\bibitem[Giakoumakis and Vanherpe, 1997]{GV-1997}
Giakoumakis, V. and Vanherpe, J.-M. (1997).
\newblock On extended {{\(P_4\)}}-reducible and extended {{\(P_4\)}}-sparse graphs.
\newblock {\em Theor. Comput. Sci.}, 180(1-2):269--286.

\bibitem[Harary and Haynes, 2000]{HH-2000}
Harary, F. and Haynes, T.~W. (2000).
\newblock Double domination in graphs.
\newblock {\em Ars Combinatoria}, 55:201--213.

\bibitem[Haynes and Hedetniemi, 2021]{HH-2021}
Haynes, T.~W. and Hedetniemi, S.~T. (2021).
\newblock Vertex sequences in graphs.
\newblock {\em Discrete Mathematics Letters}, 6:19--31.

\bibitem[Haynes et~al., 2023]{HHH-2023}
Haynes, T.~W., Hedetniemi, S.~T., and Henning, M.~A. (2023).
\newblock {\em Domination in graphs: {C}ore concepts}.
\newblock Springer.

\bibitem[Haynes et~al., 2021]{HHH-2021}
Haynes, T.~W., Hedetniemi, S.~T., Henning, M.~A., et~al. (2021).
\newblock {\em Structures of domination in graphs}, volume~66.
\newblock Springer.

\bibitem[Heggernes and Kratsch, 2007]{threshold}
Heggernes, P. and Kratsch, D. (2007).
\newblock Linear-time certifying recognition algorithms and forbidden induced subgraphs.
\newblock {\em Nordic J. of Computing}, 14(1):87–108.

\bibitem[Jamison and Olariu, 1989]{JO3}
Jamison, B. and Olariu, S. (1989).
\newblock A new class of brittle graphs.
\newblock {\em Studies in Applied Mathematics}, 81:89--92.

\bibitem[Jamison and Olariu, 1991]{JO4}
Jamison, B. and Olariu, S. (1991).
\newblock On a unique tree representation for $p_4$-extendible graphs.
\newblock {\em Discrete Applied Mathematics}, 34:151--164.

\bibitem[Pandey and Sharma, 2025]{PS-2024}
Pandey, A. and Sharma, G. (2025).
\newblock Double dominating sequences in bipartite and co-bipartite graph.
\newblock {\em Discussiones Mathematicae Graph Theory}, 45(2):545--564.

\end{thebibliography}
%\printbibliography
\end{document}